\documentclass[11pt]{amsart}
\usepackage[T1]{fontenc}
\usepackage{amsmath,amssymb,amsthm,mathtools}
\usepackage[a4paper,margin=30mm]{geometry}
\usepackage{microtype}
\usepackage[hidelinks]{hyperref}
\hypersetup{pdftitle={Zero-energy asymptotic motions near non-minimum critical points of analytic potentials},pdfauthor={}}
\newcommand{\R}{\mathbb R}
\newcommand{\Q}{\mathbb Q}
\newcommand{\N}{\mathbb N}
\newcommand{\Z}{\mathbb Z}
\newcommand{\HH}{\mathcal H}
\DeclareMathOperator{\ord}{ord}

\DeclareMathOperator{\supp}{supp}
\DeclareMathOperator{\diag}{diag}
\newtheorem{theorem}{Theorem}[section]
\newtheorem{proposition}[theorem]{Proposition}
\newtheorem{lemma}[theorem]{Lemma}

\newtheorem*{externaltheorem}{External theorem}
\theoremstyle{definition}
\newtheorem{definition}[theorem]{Definition}
\theoremstyle{remark}
\newtheorem{remark}[theorem]{Remark}
\numberwithin{equation}{section}
\allowdisplaybreaks[2]
\title{Instability of non-minimal critical points of Newton's equation with real analytic potentials}
 \author{Xiaoping Yuan}
 \address{School of Mathematical Sciences, Fudan university, Shanghai, China}
 \email{xpyuan@fudan.edu.cn}
\date{}
\subjclass[2020]{Primary 34D20; Secondary 37J25, 70H14}
\keywords{Analytic potential, Newton equation, asymptotic motion, Lyapunov instability, Hahn series, formal invariant curve}
\begin{document}
\begin{abstract}
We prove that any non-minimal critical point of a real analytic potential on a finite-dimensional Euclidean space possesses a nonconstant zero-energy Newton trajectory that converges to the equilibrium as \(t\to+\infty\). By time reversibility, this yields forward Lyapunov instability. We do not assume the critical point is isolated.

The proof constructs angular minima via quadratic regularization, matching the regularization coefficient to the leading inertial term. A free forcing coefficient allows continuation past bounded accumulations of correction orders. For polynomial potentials, this produces a complete Hahn-series solution, and differential algebraicity gives a common Puiseux denominator. A fixed-denominator continuation estimate together with a finite Taylor transfer extends the construction to analytic potentials. The realization theorem for the resulting invariant formal phase curve then yields a genuine trajectory. We also exhibit an isolated polynomial singularity for which an exact Euler-field identity admits no relatively differentiable solution on the negative set of the potential.

\end{abstract}
\maketitle

\section{Introduction}
We resolve a classical problem of Arnold concerning the instability of non-minimal critical points of Newton’s equation with real analytic potentials.
Let \(U:(\R^n,0)\to(\R,0)\) be real analytic, with \(dU(0)=0\),
and consider  Newton's equation
\begin{equation}\label{eq:newton}
 \ddot q=-\nabla U(q).
\end{equation}
The Lagrange--Dirichlet theorem gives Lyapunov stability at a strict
local minimum of the potential. Its converse includes the question
whether a critical point which is not a local minimum must be unstable.
For analytic potentials this is a classical problem; in the isolated
case it is Arnold's Problem~1971-4 \cite[p.~10]{Arnold2004} :
{\it 
Prove the instability of the equilibrium $0$ of an analytic system \(\ddot x =-\partial U/\partial x\)  in the case where the isolated (in \(\mathbb C^n\)?) critical point $0$ of the potential $U$
is not a minimum.}
The question is nontrivial even for the Euclidean kinetic energy
\(|\dot q|^2/2\). In this paper, we prove the following theorem.

\begin{theorem}\label{thm:main}
Let \(n\ge1\), and let \(U:(\R^n,0)\to(\R,0)\) be real analytic with
\(dU(0)=0\). If \(0\) is not a  minimum of \(U\), there is a
nonconstant solution \(q:[0,+\infty)\to\R^n\) of \eqref{eq:newton},
contained in a sufficiently small neighborhood of zero, such that
\begin{equation}\label{eq:mainorbit}
 (q(t),\dot q(t))\longrightarrow(0,0)\quad(t\to+\infty),\qquad
 \frac12|\dot q(t)|^2+U(q(t))=0.
\end{equation}
It may be chosen with \(U(q(t))<0\) and \(\dot q(t)\ne0\) for all
\(t\ge0\). Moreover, the equilibrium \((0,0)\) is Lyapunov unstable in forward time.
\end{theorem}

Non-minimality means that every neighborhood of zero contains a point
where \(U<0\).  The negative set can approach the equilibrium through
curved regions which are invisible to the leading homogeneous term
along fixed rays. Neither real nor complex isolation is assumed. The kinetic
energy remains Euclidean throughout; a nonlinear change of configuration
coordinates generally changes this metric. We do not assert total
instability, or an extension to arbitrary analytic kinetic metrics.
The asymptotic trajectory has one prescribed endpoint and need not be
homoclinic or heteroclinic.

We use the standard forward-time convention. The instability conclusion
means that some \(\epsilon_0>0\) satisfies
\begin{equation}\label{eq:forward-instability}
 \forall\delta>0\quad\exists q_\delta,\ T_\delta>0:\qquad
 |(q_\delta(0),\dot q_\delta(0))|<\delta,\quad
 |(q_\delta(T_\delta),\dot q_\delta(T_\delta))|>\epsilon_0,
\end{equation}
where \(q_\delta\) solves \eqref{eq:newton} on \([0,T_\delta]\).
It follows by reversing finite segments of the incoming motion, with
the velocity sign reversed as well.

\subsection{Related results}
The classical sufficient condition originates in Lagrange's work of 1788 and
Dirichlet's proof of 1846; see the historical discussion in
\cite[\S1]{Palamodov2020}. 
Lyapunov's classical criteria cover a negative direction in the
quadratic part and a negative definite first nonzero homogeneous term
\cite{Lyapunov1897}; his force function has the opposite sign from our
potential. Kozlov and Palamodov \cite{KozlovPalamodov1982} obtain an
asymptotic motion when the first nonzero homogeneous term itself has no
local minimum; see also \cite[Theorem~1]{Palamodov2020}. Such a condition
does not follow from non-minimality of the full germ. For example,
\((x-y^2)^2-y^6\) has nonnegative quadratic part, but is negative along
\((x,y)=(s^2,s)\).

The planar analytic problem was treated by Palamodov
\cite{Palamodov1977} and Taliaferro \cite{Taliaferro1980}; Kozlov
\cite{Kozlov1981} extended the isolated planar result to analytic positive
definite kinetic energy. Brunella \cite{Brunella1998} gives a geometric
approach through vector fields on invariant three-dimensional analytic
varieties, with a corollary for planar mechanics. The dimension in that
paper's title is not a statement about arbitrary three-dimensional
configuration spaces.

Higher-dimensional sufficient conditions include the finite-jet and
Hessian-rank hypotheses of Freire, Garcia and Tal
\cite{FreireGarciaTal2008}, and the radial-derivative and jet criteria
of Garcia and Alva Morales \cite{GarciaMorales2017}. Palamodov's
total-instability criterion \cite[Theorem~2]{Palamodov2020} uses a
relative \(C^1\) field on \(\{U<0\}\cup\{0\}\), nonnegativity of its
associated quadratic form at zero, and a positive lower bound for
\(A(U)/U\). It requires an inequality, rather than the identity
\(A(U)=U\). His Theorem~3 applies when the first homogeneous term has
no nonzero real critical point. Burgos \cite{Burgos2025} establishes
generic total instability for analytic electromagnetic Lagrangian
systems under a weak-magnetism hypothesis and gives further criteria
in degenerate cases. His 2025 author version describes the general
analytic non-minimum problem above two dimensions as open. Theorem~
\ref{thm:main} concerns the Euclidean, nonmagnetic system and provides
one zero-energy asymptotic motion, without these additional hypotheses.

The non-minimum question is distinct from asking whether stability
forces a \emph{strict} minimum. The latter also concerns non-strict
minima; see, for example, \cite{BurgosMadernaPaternain2021}. For
\(U(x,y)=x^2\), the equilibrium is unstable because small velocities in
the free \(y\)-direction cause escape, whereas every zero-energy motion
is constant. Thus the hypothesis in Theorem~\ref{thm:main} cannot be
replaced by failure of a strict minimum. Analyticity also matters:
smooth counterexamples to the non-minimum converse are discussed in
\cite{BurgosMadernaPaternain2021,Palamodov2020}.

\subsection{Construction and proof structure}
Write
\[
 \mu(r)=\min_{|\theta|=1}U(r\theta)
       =-\alpha r^p+o(r^p),\qquad \alpha>0.
\]
Subanalytic selection gives a convergent Puiseux curve of spherical
minima. In angular coordinates the Newton residual is the sum of the
angular potential gradient and an inertial term. We replace the
spherical minimum problem by
\[
 J(r,u)=V(r,u)+\frac\lambda2r^p|u|^2+\ell(r)\cdot u.
\]
Matching \(\lambda\) to the order of a minimizing displacement cancels
the leading residual. Global minimum comparisons provide the estimates
needed both between consecutive corrections and between arbitrary
stages.

For polynomial potentials, these minimum problems are continued in the
real closed Hahn field \(\R((r^\Q))\). At a bounded accumulation of
orders, a free coefficient of the forcing cancels the residual at the
limiting order. Countable ordinal continuation then gives a complete
formal solution. Differential algebraicity and Matusinski's support
theorem \cite{Matusinski2011}, following the differential
Newton--Puiseux work of Cano \cite{Cano1993}, imply a finite Puiseux
denominator. This support result is used only after the full equation
has been solved.

For analytic potentials, a minimum with residual order
\(f>(3p-m)/2\), where \(m=\ord U\), admits continuation with one fixed
denominator. A high Taylor polynomial and the completed polynomial
solution supply a finite forcing which reaches that threshold. The
resulting ordinary formal curve is lifted to the original Newton
phase field. The realization theorem of Le Gal and Sanz S\'anchez
\cite{LeGalSanz2026} then gives an actual trajectory. Its time
orientation and energy are determined in the original system.

The principal task is to construct an invariant irreducible real formal
curve in phase space with the leading terms of an incoming zero-energy
motion. A Puiseux selection of spherical minima provides an initial
configuration-space curve along which
\(U=-\alpha r^p+o(r^p)\), with \(\alpha>0\). This curve supplies the
negative leading potential, but generally fails to satisfy the Newton
equation: angular minimality cancels the angular potential gradient,
whereas the motion also has an angular inertial term. Writing
\(q(r)=r\phi(u(r))\), the zero-energy identity determines the radial
speed and leaves the shape equation
\[
 R[u]=\partial_uV(r,u)+\mathcal I[u]=0,
 \qquad V(r,u)=U(r\phi(u)).
\]
The word \emph{shape} means a curve of directions expressed as a
function of radius.

We construct a formal solution of this equation by successive
finite-dimensional minimum problems at each radius. The auxiliary
objectives have the form
\[
 J(r,z)=V(r,z)+\frac{\lambda}{2}r^p|z|^2+\ell(r)\cdot z.
\]
Global minimum comparisons control both the displacement of the
minimizing center and the change in potential. The linear forcing
incorporates the current Newton residual, while the quadratic
coefficient is matched to the leading inertial response of the
correction. The resulting cancellation raises the residual order
without changing the negative leading potential.

Completing these corrections requires more than a strictly increasing
sequence of orders. For polynomial potentials, we work in a rational
Hahn field to pass through possible bounded accumulations; a free
forcing coefficient cancels the residual at the limiting order.
Only after obtaining a complete solution do we use differential
algebraicity and the support theorem to obtain a common Puiseux
denominator. A finite Taylor transfer and a fixed-denominator
continuation argument treat analytic potentials. For the resulting
formal shape, define
\[
 q(r)=r\phi(\widehat u(r)),\qquad
 v(r)=-q_r(r)\sqrt{\frac{-2U(q(r))}{|q_r(r)|^2}}.
\]
After a finite ramification \(r=s^N\), the formal phase
parametrization \(\Gamma(s)=(q(s^N),v(s^N))\) satisfies
\(X_U\circ\Gamma=a(s)\Gamma'(s)\). Its underlying irreducible formal
curve is invariant and is not contained in the formal singular
locus. Theorem~\ref{thm:realization-input} then produces an actual
trajectory on the specified half-branch. The signs of its leading
position and velocity determine forward approach, and local
uniqueness excludes arrival at the equilibrium in finite time.

Sections~\ref{sec:angular}--\ref{sec:matched} establish the angular
estimates and matching. Sections~\ref{sec:hahn}--\ref{sec:puiseux}
treat polynomial potentials, and Sections~\ref{sec:threshold}--
\ref{sec:analytic-transfer} transfer the construction to analytic
potentials. Section~\ref{sec:realization} proves
Theorem~\ref{thm:main}. Section~\ref{sec:examples} checks the threshold
on curved valleys. Section~\ref{sec:euler} gives an independent
obstruction to exact Euler-field certificates; its relation to
Saito's Jacobian-ideal criterion \cite{Saito1971} is kept separate
from the trajectory construction.

\section{Angular equations and valuation estimates}\label{sec:angular}

\subsection{Selection of spherical minima}
We recall the rectilinearization input before deriving the selection
statement. Coordinate quadrants allow zero coordinates as well as
strictly positive and strictly negative coordinates.

\begin{externaltheorem}[Bierstone--Milman, {\cite[Theorem~0.2]{BM1988}}]
Let \(M\) be a real analytic manifold of pure dimension \(d\),
\(S\subset M\) subanalytic, and \(K\subset M\) compact. There are
finitely many analytic maps \(\pi_i:\R^d\to M\) and compact sets
\(L_i\subset\R^d\) such that \(\bigcup_i\pi_i(L_i)\) contains a
neighborhood of \(K\), and each \(\pi_i^{-1}(S)\) is a union of
coordinate quadrants.
\end{externaltheorem}

\begin{lemma}[Analytic curve selection]\label{lem:curve-selection}
If \(S\subset M\) is subanalytic and
\(x_0\in\overline S\setminus S\), there is an analytic map
\(\gamma:(-\epsilon,\epsilon)\to M\) with
\(\gamma(0)=x_0\) and \(\gamma(s)\in S\) for \(0<s<\epsilon\).
\end{lemma}
\begin{proof}
Apply rectilinearization with \(K=\{x_0\}\). Choose
\(x_j\in S\) tending to \(x_0\), and lifts in the compact sets
\(L_i\). After subsequences, one map \(\pi\) and one quadrant
\(Q\subset\pi^{-1}(S)\) contain lifts \(b_j\to b\in\overline Q\),
with \(\pi(b)=x_0\). For any \(a\in Q\), the segment
\(b+s(a-b)\) lies in \(Q\) for \(0<s<1\). Its image under
\(\pi\), restricted near zero, is the required analytic arc.
It is nonconstant because \(x_0\notin S\).
\end{proof}

Apply the lemma to the semianalytic negative set. The resulting arc
has the expansion
\[
 q_*(s)=s^a(v+O(s)),\qquad a\in\N_{>0},\quad v\ne0.
\]
For \(s>0\) its radius is
\[
 |q_*(s)|=c s^a h(s),\qquad c=|v|>0,\quad h(0)=1,
\]
where \(h\) is analytic and positive. This follows by taking the
positive analytic square root of
\(|s^{-a}q_*(s)|^2/c^2\). Its derivative is positive for small
positive \(s\), so the radii cover an interval \((0,r_0)\).
Consequently
\begin{equation}\label{eq:minimum}
 \mu(r):=\min_{|\theta|=1}U(r\theta)<0\qquad(0<r<r_0).
\end{equation}

The graph of \(\mu\) is subanalytic. Indeed, on bounded radius and
value intervals, both
\[
 E=\{(r,y):\exists\theta\in S^{n-1},\ U(r\theta)=y\},
 \quad
 B=\{(r,y):\exists\theta\in S^{n-1},\ U(r\theta)<y\}
\]
are projections of relatively compact semianalytic sets. Their
difference \(E\setminus B\) is precisely the minimum graph, by
compactness of the sphere. The subanalytic complement theorem
\cite[Theorem~3.10]{BM1988} applies. This is also the compact
minimization argument in \cite[Remark~3.11(2)]{BM1988}.

The set of minimizers \(\{(r,\theta):r>0, U(r\theta)=\mu(r)\}\)
is therefore subanalytic. Select a boundary point \((0,e)\) and
apply Lemma~\ref{lem:curve-selection} to obtain an analytic arc
\((\rho(s),\Theta(s))\), with \(\rho(s)>0\) for \(s>0\).
Write \(\rho(s)=c_0s^b h_0(s)\), with \(c_0>0\), \(h_0(0)=1\).
The analytic parameter change
\(z=c_0^{1/b}s h_0(s)^{1/b}\) has nonzero derivative and gives
\(\rho=z^b\). Thus \(\theta_*(r)\) is convergent in \(r^{1/b}\),
tends to \(e\), and minimizes on every sufficiently small sphere.
Since \(\mu\) is negative and convergent in this parameter,
\begin{equation}\label{eq:minimum-leading}
 \mu(r)=-\alpha r^p+o(r^p),\qquad \alpha>0,\quad p\in\Q.
\end{equation}
If \(m=\ord U\), then \(|U(r\theta)|\le Cr^m\) uniformly on the
sphere, whence \(p\ge m\ge2\).

The same argument will be used for auxiliary minimum problems.
\begin{lemma}[Compact minimum selection]\label{lem:compact-selection}
Let \(K\) be a nonempty compact semianalytic set, and let \(F(r,z)\) be
continuous in \(z\in K\) for \(r>0\). Suppose that after a finite
ramification of \(r\) and multiplication by a positive monomial,
\(F\) extends analytically to a neighborhood of \(\{0\}\times K\).
Then its global minima on \(K\) admit a convergent real Puiseux
selection for all sufficiently small positive \(r\).
\end{lemma}
\begin{proof}
The multiplier does not change the minimizers. The preceding compact
projection and complement argument makes their set subanalytic.
Compactness supplies a boundary point at radius zero; analytic curve
selection and the root reparametrization just used give the result.
\end{proof}

\subsection{Radial solutions and the angular equation}\label{subsec:radial}
If a fixed unit vector \(e\) is angularly critical for every small
\(r>0\), with \(U(re)<0\), solve
\[
 \dot r=-\sqrt{-2U(re)},\qquad q=re.
\]
Then \(\ddot r=-\nabla U(re)\cdot e\), while angular criticality
makes \(\nabla U(re)\) parallel to \(e\). Hence \(q\) solves the
full Newton equation and has zero energy. Its leading potential is
\(-\alpha_e r^{p_e}\), with \(p_e\ge2\), and
\(\int_0^{r_0}(-2U(re))^{-1/2}\,dr=+\infty\). The motion tends
to the equilibrium in positive infinite time. This covers constant
minimum directions, all later radial alternatives, and the case
\(n=1\). In general \(p_e\) need not equal \(p\).

Assume henceforth \(n\ge2\), put \(k=n-1\), and rotate \(e\) to
\(e_n\). Use the rational chart
\begin{equation}\label{eq:chart}
 \rho(u)=1+\frac{|u|^2}{4},\qquad
 \phi(u)=\frac{(u,1-|u|^2/4)}{\rho(u)},\qquad g(u)=\rho(u)^{-2}.
\end{equation}
Then \(|\phi|=1\), \(\phi_u^T\phi_u=gI_k\), and
\(|dq|^2=dr^2+r^2g(u)|du|^2\) for \(q=r\phi(u)\).
Set
\[
 V(r,u)=U(r\phi(u))=r^mW(r,u),\qquad D=r\frac{d}{dr}.
\]
Every fixed angular derivative of \(V\) is \(O(r^m)\) uniformly on
a sufficiently small fixed closed angular ball \(K\). The selected
minimum curve lies in its interior for small \(r\). The Lagrangian in these coordinates is
\[
 L=\frac12\dot r^{\,2}+\frac12r^2g(u)|\dot u|^2-V(r,u).
\]
Its angular Euler--Lagrange equations are
\[
 \frac{d}{dt}(r^2g\dot u_a)
 -\frac{r^2}{2}\partial_a g\,|\dot u|^2+\partial_aV=0.
\]

For a shape \(u(r)\), define the functions of radius
\begin{equation}\label{eq:energy}
 \begin{aligned}
 Q[u](r)&:=1+g(u(r))|Du(r)|^2,\\
 w[u](r)&:=-\frac{2V(r,u(r))}{Q[u](r)}.
 \end{aligned}
\end{equation}
For an actual zero-energy motion with this shape,
\(\dot r(t)^2=w[u](r(t))\); on an incoming segment with \(V<0\),
\(\dot r(t)=-\sqrt{w[u](r(t))}\). The functional \(w[u]\) is
not itself a time derivative.

The first-kind Christoffel tensor of the angular metric is
\[
 \Gamma_{\rm cov}(u)[y,y]
 =\rho(u)^{-3}\left(\frac12u|y|^2-(u\cdot y)y\right).
\]
Define
\begin{align}
 \mathcal I[u]&=w[u]g(u)D^2u
 +\left(w[u]+\frac12Dw[u]\right)g(u)Du
 +w[u]\Gamma_{\rm cov}(u)[Du,Du],\label{eq:inertia}\\
 R[u]&=\partial_uV(r,u)+\mathcal I[u].\label{eq:residual}
\end{align}
To derive the equation, substitute \(\dot u=(\dot r/r)Du\) into
\[
 \frac{d}{dt}(r^2g\dot u_a)
 -\frac{r^2}{2}\partial_a g\,|\dot u|^2+\partial_aV=0.
\]
The identities \(\ddot r=w_r/2\) and \(r\ddot r=Dw/2\) give
\(R[u]=0\). Both signs of \(\dot r\) give the same shape equation.
For \(u=a r^\sigma+o(r^\sigma)\) and
\(V(r,u)=-\alpha r^p+o(r^p)\), its leading inertia is
\(C(\sigma)a r^{p+\sigma}\), where
\begin{equation}\label{eq:C}
 C(\sigma)=2\alpha\sigma\left(\sigma+1+\frac p2\right).
\end{equation}
This function is positive and strictly increasing for \(\sigma>0\).

\subsection{Series and inertia comparison}
Let \(\HH_\Q=\R((r^\Q))\) be the ordered Hahn field, with
well-ordered rational supports and valuation \(\nu\), and put
\(K_N=\operatorname{Frac}\R\{r^{1/N}\}\). Convergent Puiseux germs
are identified with their Hahn expansions. Ordinary formal Puiseux
series belong to \(\bigcup_N\R((r^{1/N}))\); they need not converge.
For vectors, \(\nu\) is the least component valuation, equal to that
of their Euclidean norm. We set \(\nu(0)=+\infty\).

The Euler derivation acts termwise by \(D(r^q)=q r^q\), preserves
summable Hahn families, and satisfies \(\nu(Dh)\ge\nu(h)\), with
equality at positive valuation. We write \(O(r^a)\) for valuation
at least \(a\), and \(o(r^a)\) for valuation greater than \(a\),
even when \(a\in\R\setminus\Q\). The coefficient \([h]_a\) is
used only for rational \(a\). Analytic substitution at positive-order
arguments is well defined; in the polynomial Hahn construction the
chart expressions are rational. Ordinary formal substitution is used
after obtaining a common denominator.

\begin{lemma}\label{lem:inertia}
Let \(u,v\) be convergent Puiseux or Hahn shapes for a fixed chart
potential. Suppose \(\nu(u)=\nu(v)=\tau>0\),
\(V(r,u)=-\alpha r^p+o(r^p)\), and, for some real \(s>\tau\),
\[
 h=v-u=O(r^s),\qquad V(r,v)-V(r,u)=O(r^{p+\tau+s}).
\]
Then
\begin{equation}\label{eq:inertia-difference}
 \mathcal I[v]-\mathcal I[u]
 =2\alpha r^p\left(D^2h+\left(1+\frac p2\right)Dh\right)
       +o(r^{p+s}).
\end{equation}
In particular the difference has valuation at least \(p+s\); if
\(s\in\Q\), its coefficient at \(p+s\) is \(C(s)[h]_s\).
\end{lemma}
\begin{proof}
Since \(g=1+O(|u|^2)\),
\(g(v)-g(u)=O(r^{\tau+s})\) and
\(Q[v]-Q[u]=O(r^{\tau+s})\). Equation~\eqref{eq:energy} gives
\[
 w[v]-w[u]
 =-\frac{2(V(r,v)-V(r,u))}{Q[v]}
  +\frac{2V(r,u)(Q[v]-Q[u])}{Q[v]Q[u]}
 =O(r^{p+\tau+s}).
\]
The same bound holds after applying \(D\). In the first inertia term,
split the difference as
\[
 w[u]g(u)D^2h+
 (w[v]g(v)-w[u]g(u))D^2v.
\]
The second summand is \(O(r^{p+s+2\tau})\); the first has leading
part \(2\alpha r^pD^2h\). The middle inertia term similarly has
leading difference \(2\alpha(1+p/2)r^pDh\), with coefficient-change
error \(O(r^{p+s+2\tau})\). Finally
\(\Gamma_{\rm cov}(u)[Du,Du]=O(r^{3\tau})\), and its difference
at \(u,v\) is \(O(r^{s+2\tau})\). The connection contribution is
therefore also \(O(r^{p+s+2\tau})\). The remaining leading-coefficient
errors are \(o(r^{p+s})\), proving the assertion.
\end{proof}

\begin{remark}\label{rem:inertia-size}
If \(u=a r^\tau+o(r^\tau)\), then
\([\mathcal I[u]]_{p+\tau}=C(\tau)a\), so
\(\nu(\mathcal I[u])=p+\tau\). The stronger potential-value estimate
in Lemma~\ref{lem:inertia} is essential: position closeness alone
only gives the generic Taylor bound \(O(r^{m+s})\).
\end{remark}

\section{Matched minima and finite corrections}\label{sec:matched}

Throughout this section \(K\) is a fixed closed Euclidean ball about
\(0\in\R^k\), contained in the angular chart.  The spherical minimizing
curve lies in its interior for sufficiently small \(r>0\).
The estimates for the angular potential on this ball are
\[
 \|\partial_u V(r,\cdot)\|+\|\partial_u^2V(r,\cdot)\|=O(r^m).
\]
Each higher fixed angular derivative is also \(O(r^m)\).
All assertions about convergent germs below are understood for sufficiently
small positive \(r\); this interval may depend on the finitely many
parameters being compared.

We first record the elementary fixed point observation used to choose
the regularization coefficients.
\begin{lemma}\label{lem:monotone-fp}
Every nondecreasing map \(T:[a,b]\to[a,b]\) has a fixed point.
\end{lemma}
\begin{proof}
The set \(S=\{x:T(x)\le x\}\) contains \(b\).  Let \(c=\inf S\).
For \(x\in S\), monotonicity gives \(T(c)\le T(x)\le x\), whence
\(T(c)\le c\).  If \(T(c)<c\), choose \(y\) strictly between them.
Then \(T(y)\le T(c)<y\), so \(y\in S\), contradicting the definition
of \(c\).  The degenerate interval case is immediate.
\end{proof}

\begin{proposition}[Initial matching]\label{prop:initial}
Either the Newton system has an exact zero-energy radial shape, or
there exist a nonzero convergent real Puiseux germ \(u_0\), a number
\(\tau=\nu(u_0)>0\), and \(\lambda_0=C(\tau)\), such that \(u_0\)
globally minimizes
\[
 J_0(r,z)=V(r,z)+\frac{\lambda_0}{2}r^p|z|^2
 \quad\text{on }K,
\]
and
\[
 V(r,u_0)=-\alpha r^p+o(r^p),\qquad
 \nu (R[u_0])>p+\tau .
\]
The last inequality includes the possibility \(R[u_0]=0\).
\end{proposition}
\begin{proof}
Let \(u_*(r)\to0\) be the angular representation of the spherical
minimizing curve.  A constant minimizing direction gives an exact
radial shape, so suppose \(u_*\not\equiv0\).
For each real \(\lambda>0\), choose a convergent Puiseux global
minimizer \(u_\lambda\) of
\[
 J_\lambda(r,z)=V(r,z)+\frac{\lambda}{2}r^p|z|^2.
\]
Such choices follow from Lemma~\ref{lem:compact-selection}.
Since \(V(r,z)\ge\mu(r)=V(r,u_*)\), comparison with \(u_*\) yields
\begin{equation}\label{eq:initial-comparison}
 |u_\lambda|\le |u_*|,\qquad
 \mu(r)\le V(r,u_\lambda)
 \le\mu(r)+\frac{\lambda}{2}r^p|u_*|^2 .
\end{equation}
Consequently \(u_\lambda\to0\), it is an interior minimizer, and its
potential has the required leading term.  If one of these germs is
zero, its stationary equation gives \(\partial_uV(r,0)=0\), while
\eqref{eq:initial-comparison} gives \(V(r,0)<0\).  The axis is then
an exact radial shape.

Otherwise set \(\delta(\lambda)=\nu(u_\lambda)>0\).
Adding the two minimizing inequalities for \(\lambda_2>\lambda_1\)
shows that \(|u_{\lambda_2}|\le|u_{\lambda_1}|\).
Thus \(\delta(\lambda)\) is nondecreasing, independently of
nonuniqueness in the choices.  Equation \eqref{eq:initial-comparison}
gives the positive uniform lower bound
\(\delta(\lambda)\ge\nu(u_*)\).

There is also a finite upper bound independent of \(\lambda\).
Put \(g_0(r)=\partial_uV(r,0)\).  If \(g_0\ne0\), the stationary
equation at \(u_\lambda\), together with the Hessian estimate, gives
\[
 |g_0(r)|\le C_\lambda r^m|u_\lambda(r)|,
 \qquad \delta(\lambda)\le\nu(g_0)-m .
\]
If \(g_0=0\) and \(V(r,0)<0\), the axis again solves the problem.
In the remaining case \(V(r,0)\ge0\), and
\[
 c_0r^p\le V(r,0)-V(r,u_\lambda)
 \le C r^m|u_\lambda|
\]
for a fixed \(c_0>0\), so \(\delta(\lambda)\le p-m\).
We have therefore obtained
\(0<\underline\delta\le\delta(\lambda)\le\overline\delta<\infty\).

The nondecreasing map
\(\lambda\mapsto C(\delta(\lambda))\) sends
\([C(\underline\delta),C(\overline\delta)]\) into itself.
Lemma~\ref{lem:monotone-fp} provides a fixed point \(\lambda_0\).
Let \(u_0=u_{\lambda_0}\) and \(\tau=\delta(\lambda_0)\).
Stationarity gives
\[
 R[u_0]=\mathcal I[u_0]-C(\tau)r^pu_0 .
\]
The leading-term inertia calculation in Lemma~\ref{lem:inertia}
shows that this expression has valuation strictly greater than
\(p+\tau\).
\end{proof}

\begin{definition}\label{def:matched}
Fix the leading data \(\alpha,p,\tau\), with \(\tau>0\).
A \emph{matched minimum} is a triple \((u,d,\ell)\) such that
\begin{enumerate}
\item \(u\) globally minimizes
\[
 J(r,z)=V(r,z)+\frac{\lambda}{2}r^p|z|^2+\ell(r)\cdot z
 \quad\text{on }K,\qquad \lambda=C(d);
\]
\item \(\nu(u)=\tau\), \(d\ge\tau\), and
 \(V(r,u)=-\alpha r^p+o(r^p)\);
\item \(f:=\nu( R[u])>p+d\).
\end{enumerate}
The label \(d\) is real and need not equal the valuation of the
center \(u\).  An exact center, for which \(f=+\infty\), is terminal.
For convergent data the minimum is interpreted for real \(r>0\);
for polynomial potentials it will also be interpreted over the
ordered Hahn field introduced below.
\end{definition}

\begin{proposition}[Successor correction]\label{prop:successor}
Let \((u,d,\ell)\) be a matched minimum with finite residual
valuation \(f\).  There is another matched minimum
\((u^+,d^+,\ell^+)\) such that
\begin{align}
 d^+&=\nu(u^+-u),&
 f-p&\le d^+\le f-m, \label{eq:successor-orders}\\
 \lambda^+&=C(d^+)>\lambda,&
 \ell^+&=\ell+R[u]-(\lambda^+-\lambda)r^pu. \label{eq:successor-force}
\end{align}
In particular \(d^+>d\) and
\begin{equation}\label{eq:successor-improvement}
 \nu( R[u^+])>p+d^+\ge f.
\end{equation}
For convergent Puiseux data, the new data are convergent Puiseux
germs, with a possibly larger common denominator.
The same assertion holds in the Hahn setting for polynomial \(U\).
\end{proposition}
\begin{proof}
For a fixed real \(c>0\), choose a global minimizer \(u+h_c\) of
\[
 J(r,z)+\frac c2r^p|z-u|^2+R[u]\cdot(z-u).
\]
Existence in the convergent case follows as in
Proposition~\ref{prop:initial}; the Hahn version follows from
the ordered-field facts in Section~\ref{sec:hahn}.
Comparison with \(u\), and the global minimality of \(u\) for \(J\),
give
\begin{equation}\label{eq:correction-bound}
 |h_c|\le\frac{2|R[u]|}{cr^p}.
\end{equation}
Thus \(\sigma(c):=\nu(h_c)\ge f-p>d\ge\tau\); in particular
the new minimum is interior.  The displacement is nonzero because
the gradient of the new objective at \(u\) is \(R[u]\ne0\).
The Hessian of this objective is \(O(r^m)\), so comparison of its
gradients at \(u\) and \(u+h_c\) gives
\[
 |R[u]|\le C_c r^m|h_c|,\qquad \sigma(c)\le f-m.
\]
Adding the minimizing inequalities at two parameters
\(c_2>c_1\) gives \(|h_{c_2}|\le|h_{c_1}|\), for any choices
of minimizers.  Hence \(\sigma(c)\) is nondecreasing.
The map \(c\mapsto C(\sigma(c))-\lambda\) sends the positive
closed interval
\[
 [C(f-p)-\lambda,\ C(f-m)-\lambda]
\]
into itself.  Choose a fixed point \(c\) by
Lemma~\ref{lem:monotone-fp}, and put
\[
 h=h_c,\quad u^+=u+h,\quad d^+=\sigma(c),\quad
 \lambda^+=\lambda+c=C(d^+).
\]
After discarding a constant independent of \(z\), the new objective
has forcing \(\ell^+\) in \eqref{eq:successor-force}.

Write \(\Delta J=J(r,u+h)-J(r,u)\).  The minimizing comparison gives
\[
 0\le\Delta J\le |R[u]|\,|h|.
\]
Since \(\partial_uV(r,u)=R[u]-\mathcal I[u]\) has valuation at
least \(p+\tau\), the exact quadratic expansion
\begin{equation}\label{eq:energy-expansion-successor}
 V(r,u+h)-V(r,u)
 =\Delta J+\partial_uV(r,u)\cdot h-\frac{\lambda}{2}r^p|h|^2
\end{equation}
shows that
\[
 \nu\bigl(V(r,u+h)-V(r,u)\bigr)\ge p+\tau+\sigma(c).
\]
It follows that the center and energy retain their leading terms.
Subtracting stationary equations gives the exact identity
\begin{equation}\label{eq:exact-successor-residual}
 R[u+h]=\mathcal I[u+h]-\mathcal I[u]-\lambda^+r^ph.
\end{equation}
If \(h=a r^{d^+}+o(r^{d^+})\), Lemma~\ref{lem:inertia}
identifies the leading inertia difference as
\(C(d^+)a r^{p+d^+}\).  Our choice
\(\lambda^+=C(d^+)\) cancels it, proving
\eqref{eq:successor-improvement}.
\end{proof}

\section{Hahn minima and continuation through accumulation}\label{sec:hahn}

For polynomial \(U\), the chart potential \(V(r,z)\) is rational
in \(r,z\), with denominators positive on \(K\).  We work in
\[
 \HH=\R((r^{\Q})).
\]
An element of this field has well-ordered support in \(\Q\), with
no local-finiteness requirement.  The field is ordered by its
leading coefficient and is real closed: its coefficient field is
real closed and its value group is divisible
\cite[Section~1.2]{DMM1994}.  The termwise Euler derivation
\(D(r^q)=q r^q\) is a strong derivation.  For nonzero positive
valuation it preserves valuation; in general it does not decrease
valuation.

We use two elementary consequences of real-closed-field quantifier
elimination.  First, a continuous rational function whose
denominator has no zero on the closed bounded semialgebraic ball
\(K_{\HH}\) attains its minimum there.  The coefficients of the
function may be arbitrary elements of \(\HH\).
This is transfer of the semialgebraic extreme-value statement,
not compactness of \(K_{\HH}\) in its order topology.
Interior minimizers satisfy the algebraic stationary equations.
The uniform rational derivative bounds used above transfer as well.

Second, a convergent Puiseux minimizer of such a rational objective
remains a global minimizer against every Hahn test point.
Indeed, its coefficients, the chosen center, and \(r^p\) enter the
global minimum assertion as finitely many field parameters.
Quantifier elimination reduces this assertion to finitely many
polynomial equalities and sign conditions on these parameters.
Each resulting convergent Puiseux expression has the same eventual
sign on the positive real axis as its Hahn expansion.
This argument does not require these parameter germs to be algebraic
over \(\R(r)\).

All regularization coefficients, including \(C(d)\) for an irrational
label \(d\), are real constants.  No monomial of irrational exponent
is introduced by such a label.

\subsection{Comparison of nonconsecutive minima}

\begin{lemma}[Coherence comparison]\label{lem:coherence}
Let \((u,d,\ell)\) be a matched minimum with finite residual
valuation \(f\), and let \(v\) globally minimize an objective
\[
 J'(r,z)=V(r,z)+\frac{\lambda'}2r^p|z|^2+\ell'(r)\cdot z,
 \qquad \lambda'>\lambda=C(d).
\]
Put
\[
 E=\ell'-\ell+(\lambda'-\lambda)r^pu.
\]
If \(\nu(E)\ge f\), then \(v\) is interior and, with \(h=v-u\),
\begin{align}
 \nu(h)&\ge f-p, \label{eq:coherence-position}\\
 \nu\bigl(V(r,v)-V(r,u)\bigr)&\ge f+\tau, \label{eq:coherence-energy}\\
 \nu\bigl(\mathcal I[v]-\mathcal I[u]\bigr)&\ge f. \label{eq:coherence-inertia}
\end{align}
In fact, when \(h\ne0\), the energy bound improves to
\(p+\tau+\nu(h)\).
\end{lemma}
\begin{proof}
Up to a constant independent of \(z\), one has
\[
 J'(r,z)=J(r,z)+\frac{\lambda'-\lambda}{2}r^p|z-u|^2
       +E\cdot(z-u).
\]
Comparison at the two minimizing points gives
\[
 |h|\le\frac{2|E|}{(\lambda'-\lambda)r^p},\qquad
 0\le J(r,v)-J(r,u)\le |E|\,|h|.
\]
If \(E=0\), positive quadratic growth forces \(h=0\).
Otherwise \(\nu(h)\ge f-p>\tau\), so \(v\to0\) and is interior.
The exact expansion
\[
 V(r,v)-V(r,u)
 =J(r,v)-J(r,u)+\partial_uV(r,u)\cdot h
       -\frac{\lambda}{2}r^p|h|^2
\]
and \(\nu(\partial_uV(r,u))\ge p+\tau\) yield the stronger
energy estimate.  Lemma~\ref{lem:inertia}, with
\(s=f-p>\tau\), proves \eqref{eq:coherence-inertia}.
\end{proof}

For an ordered family of matched minima we maintain, for every
\(\beta<\gamma\), the forcing invariant
\begin{equation}\label{eq:all-pair-forcing}
 E_{\beta\gamma}
 =\ell_\gamma-\ell_\beta+
   (\lambda_\gamma-\lambda_\beta)r^pu_\beta,\qquad
 \nu(E_{\beta\gamma})\ge f_\beta.
\end{equation}
Lemma~\ref{lem:coherence} then supplies all three comparison
estimates with the corresponding indices.
At a successor, \(E_{\gamma,\gamma+1}=R[u_\gamma]\).
The exact composition rule
\begin{equation}\label{eq:forcing-composition}
 E_{\beta\eta}=E_{\beta\gamma}+E_{\gamma\eta}
       +(\lambda_\eta-\lambda_\gamma)r^p(u_\beta-u_\gamma)
\end{equation}
shows that the invariant extends from all preceding pairs to a
new successor.  Indeed, residual valuations increase, and the last
term has valuation at least \(f_\beta\) by
\eqref{eq:coherence-position}.
This all-pair formulation is needed at limits; adjacent comparison
alone would not suffice.

\subsection{A bounded countable limit}

\begin{proposition}[Passage through a bounded accumulation]\label{prop:bounded-limit}
Let \(\eta\) be a countable limit ordinal, and suppose matched
minima have been constructed at all indices \(\beta<\eta\), with
finite residual valuations, strictly increasing labels, the
successor bounds \eqref{eq:successor-orders}--\eqref{eq:successor-improvement},
and the all-pair invariant \eqref{eq:all-pair-forcing}.
If \(L=\sup_{\beta<\eta}d_\beta<\infty\), there is a matched
minimum at index \(\eta\), with label \(d_\eta=L\), which preserves
the all-pair invariant and all comparison estimates.
\end{proposition}
\begin{proof}
Suppose the construction has reached every index below a countable
limit ordinal \(\eta\), without obtaining an exact center, and that
\[
 d_\beta\uparrow L<\infty\quad(\beta\to\eta).
\]
At each successor,
\[
 p+d_\beta<f_\beta\le p+d_{\beta+1}.
\]
Every \(\beta<\eta\) has its successor below \(\eta\), hence
\begin{equation}\label{eq:limit-cutoff}
 f_\beta\uparrow F:=p+L,\qquad f_\beta<F.
\end{equation}
Let
\[
 \Lambda=C(L),\qquad
 \epsilon_\beta=\Lambda-\lambda_\beta>0,\qquad
 b_\beta=\ell_\beta-\epsilon_\beta r^pu_\beta.
\]
The old center has global quadratic growth for the objective
with fixed coefficient \(\Lambda\):
\begin{align}
 &V(r,z)+\frac{\Lambda}{2}r^p|z|^2+b_\beta\cdot z
       -V(r,u_\beta)-\frac{\Lambda}{2}r^p|u_\beta|^2
       -b_\beta\cdot u_\beta \notag\\
 &\hspace{35mm}\ge
       \frac{\epsilon_\beta}{2}r^p|z-u_\beta|^2 .
 \label{eq:fixed-quadratic-growth}
\end{align}
Moreover, for \(\beta<\gamma<\eta\),
\begin{equation}\label{eq:fixed-forcing-coherence}
 b_\gamma-b_\beta
 =E_{\beta\gamma}-(\Lambda-\lambda_\gamma)r^p(u_\gamma-u_\beta),
 \qquad \nu(b_\gamma-b_\beta)\ge f_\beta.
\end{equation}

Consequently the coefficients at every rational exponent less than
\(F\) stabilize.  Set all coefficients at or above \(F\) to zero
and denote the resulting vector by \(b_{<F}\).
It is a Hahn vector.  To verify well ordering in one coordinate,
take a nonempty subset of its support and an exponent \(q\) in that
subset.  Choose \(\beta\) with \(f_\beta>q\).
Below \(f_\beta\) the stabilized support is the support of a fixed
initial segment of \(b_\beta\).  The part of the chosen subset below
that cutoff has a least element, and this element is least in the
whole subset.  This proves well ordering, without asserting local
finiteness.  Also
\begin{equation}\label{eq:prefix-error}
 \nu(b_{<F}-b_\beta)\ge f_\beta\quad(\beta<\eta).
\end{equation}

Assume first that \(L\in\Q\), so that \(F\in\Q\).  For each
\(a\in\R^k\), choose a global Hahn minimizer \(u_a\) of
\begin{equation}\label{eq:boundary-objective}
 V(r,z)+\frac{\Lambda}{2}r^p|z|^2+
                 (b_{<F}+a r^F)\cdot z .
\end{equation}
For every fixed \(\beta<\eta\), its forcing difference from the
fixed-quadratic old objective is of valuation at least \(f_\beta\).
Lemma~\ref{lem:coherence}, or its minimizing comparison with
\eqref{eq:fixed-quadratic-growth}, therefore gives
\begin{align}
 \nu(u_a-u_\beta)&\ge f_\beta-p, \label{eq:limit-position}\\
 \nu\bigl(V(r,u_a)-V(r,u_\beta)\bigr)&\ge f_\beta+\tau, \label{eq:limit-energy}\\
 \nu\bigl(\mathcal I[u_a]-\mathcal I[u_\beta]\bigr)&\ge f_\beta .
 \label{eq:limit-inertia}
\end{align}
In particular \(u_a\) is interior and retains both leading terms.
The constants \(\epsilon_\beta\) may tend to zero as
\(\beta\to\eta\); this causes no loss, since each comparison fixes
\(\beta\), and every nonzero real constant has valuation zero.

Subtracting stationary equations yields
\begin{align}
 R[u_a]-R[u_\beta]
 &=\mathcal I[u_a]-\mathcal I[u_\beta]
   -\Lambda r^p(u_a-u_\beta)
   -(b_{<F}+a r^F-b_\beta).
 \label{eq:limit-residual-comparison}
\end{align}
Every term on the right has valuation at least \(f_\beta\).
Letting \(\beta\) approach \(\eta\) shows that \(\nu (R[u_a])\ge F\).
For two choices \(a,b\), the same estimates give
\begin{equation}\label{eq:boundary-pair}
 \nu(u_a-u_b)\ge L,\qquad
 \nu\bigl(V(r,u_a)-V(r,u_b)\bigr)\ge p+\tau+L.
\end{equation}
Lemma~\ref{lem:inertia} now identifies the boundary coefficient:
\[
 [\mathcal I[u_a]-\mathcal I[u_b]]_F
       =C(L)[u_a-u_b]_L.
\]
On the other hand, subtraction of the exact stationary equations
of \eqref{eq:boundary-objective} gives
\[
 \partial_uV(r,u_a)-\partial_uV(r,u_b)
       =-\Lambda r^p(u_a-u_b)-(a-b)r^F.
\]
Since \(\Lambda=C(L)\), the unknown center coefficient cancels:
\begin{equation}\label{eq:boundary-cancellation}
 \boxed{\ [R[u_a]]_F-[R[u_b]]_F=-(a-b).\ }
\end{equation}
This identity holds for every pair of minimizing choices; it
requires neither uniqueness of a minimizer nor continuous
dependence on \(a\).
First choose a minimizer at \(a=0\), and then set
\[
 a=[R[u_{a=0}]]_F.
\]
Every corresponding minimizer has either zero residual or
residual valuation strictly greater than \(F\).
Choose one and define
\[
 u_\eta=u_a,\qquad d_\eta=L,\qquad
 \lambda_\eta=\Lambda,\qquad \ell_\eta=b_{<F}+a r^F.
\]
This is a matched minimum.  The old-pair comparisons hold by
\eqref{eq:limit-position}--\eqref{eq:limit-inertia}, and
\[
 E_{\beta\eta}=\ell_\eta-b_\beta
\]
has valuation at least \(f_\beta\) by \eqref{eq:prefix-error}.
Thus the forcing invariant also holds for every pair ending at
the new limit stage.

If \(L\notin\Q\), minimize \eqref{eq:boundary-objective} with its
last forcing term omitted, that is, with forcing \(b_{<F}\).
The preceding comparison gives \(\nu (R[u])\ge F\).
A nonzero Hahn residual has rational valuation, whereas \(F\)
is irrational, so the inequality is automatically strict.
Set \(d_\eta=L\), \(\lambda_\eta=C(L)\), and
\(\ell_\eta=b_{<F}\).  All invariants are verified as before.
The label may be irrational, but every series remains in
\(\R((r^{\Q}))\).
\end{proof}

\subsection{Completion}

\begin{proposition}[Exact polynomial Hahn shape]\label{prop:hahn}
For a polynomial potential, every matched Hahn minimum can be
continued to an exact Hahn shape \(\widehat u\) with
\[
 R[\widehat u]=0,\qquad
 \nu(\widehat u)=\tau,\qquad
 V(r,\widehat u)=-\alpha r^p+o(r^p).
\]
The construction consists of matched minima at successor and
bounded countable limit stages, and possibly a terminal unbounded
countable limit.  For every nonterminal stage \(\beta\),
\begin{align}
 \nu(\widehat u-u_\beta)&\ge f_\beta-p, \label{eq:terminal-position}\\
 \nu\bigl(V(r,\widehat u)-V(r,u_\beta)\bigr)&\ge f_\beta+\tau,
 \label{eq:terminal-energy}\\
 \nu\bigl(\mathcal I[\widehat u]-\mathcal I[u_\beta]\bigr)&\ge f_\beta .
 \label{eq:terminal-inertia}
\end{align}
\end{proposition}
\begin{proof}
An initially exact minimum needs no continuation.  Otherwise use
Proposition~\ref{prop:successor} at successors and the preceding
bounded-limit construction at every countable limit with bounded
labels.  Maintain simultaneously the strictly increasing labels,
the increasing finite residual valuations, and
\eqref{eq:all-pair-forcing}.
Equations \eqref{eq:forcing-composition} and
\eqref{eq:prefix-error} prove the successor and limit induction,
respectively.  Lemma~\ref{lem:coherence} supplies all comparisons.

Suppose a countable limit has unbounded preceding labels.
The residual valuations are then unbounded, and
\eqref{eq:coherence-position} makes the centers Cauchy for the
valuation topology: below each finite exponent cutoff, their
coefficients eventually coincide.  The stable union of these
initial segments has well-ordered support by the same least-element
argument used above.  It defines a Hahn vector \(\widehat u\),
with the original positive leading order.
The Euler derivation is valuation-nondecreasing; rational
substitution and inversion are continuous at the nonzero
denominators occurring here.  Thus \(R[u_\beta]\to0\) and
\(u_\beta\to\widehat u\) imply \(R[\widehat u]=0\).
The leading energy and the three terminal comparison estimates
pass to this limit.

If neither an exact center nor an unbounded countable limit
occurred, transfinite recursion would construct a matched minimum
at every countable ordinal, with strictly increasing real labels
\((d_\beta)_{\beta<\omega_1}\).  For each \(\beta<\omega_1\),
choose a rational number in \((d_\beta,d_{\beta+1})\).
These intervals are disjoint, so the choices would inject an
uncountable set into \(\Q\), a contradiction.
Hence one of the two terminal events occurs.
When termination occurs at a matched stage, the terminal estimates
are precisely Lemma~\ref{lem:coherence}; when it occurs at an
unbounded limit, they were proved above.
\end{proof}

\begin{remark}
This is an existence argument using countable ordinal continuation.
It neither asserts that an ordinary sequence of finite corrections
has unbounded orders, nor provides a finite algorithm.  At a bounded
accumulation it solves a new full minimum problem before
differentiating.  The boundary coefficient in
\eqref{eq:boundary-cancellation} is free because only the forcing
strictly below \(F\) was fixed by previous stages.
\end{remark}

\section{From polynomial Hahn solutions to formal Puiseux solutions}\label{sec:puiseux}

The support theorem is applied to a complete solution of the original
shape equation. An intermediate minimum still contains an auxiliary
forcing and need not have a finite Puiseux denominator.

\begin{lemma}\label{lem:scalar-support}
If \(0\ne y\in\HH_\Q\) has positive valuation and is differentially
algebraic over \(\R(r)\) for \(D=r\,d/dr\), then
\(y\in\R[[r^{1/N}]]\) for some \(N\ge1\).
\end{lemma}
\begin{proof}
Embed \(\HH_\Q\) in the rank-one field \(\HH_\R=\R((r^\R))\).
The derivation \(D(r^a)=a r^a\) is strong and satisfies the generalized
Leibniz rule. For \(\nu(f)\ne0\),
\[
 \nu(Df)=\nu(f),\qquad \nu(Df/f)=0.
\]
These verify the Hardy-type hypotheses of
\cite[Definition~2.2]{Matusinski2011}; the strict comparison between
distinct Archimedean classes is vacuous in rank one. Every \(D^j y\)
has the same positive valuation as \(y\).

Clear denominators in a nonzero scalar differential polynomial
annihilating \(y\). Its coefficient supports are finite subsets of
\(\N\), and \(\supp(Dr/r)=\{0\}\). Matusinski's Theorems~1.2
and 2.14 place the support of an exact solution in a set obtained from
these supports by finitely many elementary transformations. The
alternative of a stabilized nonzero residual in Theorem~2.14 cannot
occur for an exact solution. The transformations are sums, additive
semigroup generation, adjoining a generator, and translation of a tail
\cite[Definition~2.13]{Matusinski2011}. Each preserves containment in a
finitely generated additive subgroup after adjoining finitely many
constants. Consequently
\[
 \supp y\subset G=\Z\gamma_1+\cdots+\Z\gamma_b\subset\R.
\]
Here and below \(\Z\) denotes the integers.

The group \(G\) need not be discrete. Nevertheless, \(G\cap\Q\)
has a common denominator. If this intersection is nonzero, choose a
basis \(1,\beta_2,\ldots,\beta_d\) of the rational span of the
\(\gamma_i\), and write
\(\gamma_i=a_{i1}+\sum_{j=2}^d a_{ij}\beta_j\), with
\(a_{ij}\in\Q\). Let \(N\) clear the denominators of the
\(a_{i1}\). In a rational integer combination of the \(\gamma_i\),
all coefficients of \(\beta_j\), \(j\ge2\), vanish. Thus
\(G\cap\Q\subset(1/N)\Z\). If the intersection is zero the
assertion is immediate. Since \(y\) has rational support and positive
valuation, the result follows.
\end{proof}

\begin{proposition}\label{prop:puiseux}
Let \(U\) be polynomial. If \(u\in\HH_\Q^k\) is a positive-order
solution of \(R[u]=0\) with \(V(r,u)=-\alpha r^p+o(r^p)\), then
\(u\in\R[[r^{1/N}]]^k\) for some \(N\ge1\).
\end{proposition}
\begin{proof}
Put \(y=Du\), \(Q=1+g|y|^2\), and \(w=-2V/Q\). The potential
\(V\) is rational in \((r,u)\), since both \(U\) and the chart
\(\phi\) are rational. Account for the acceleration inside \(Dw\):
\[
 DQ=(\nabla g\cdot y)|y|^2+2g\,y\cdot Dy,\qquad
 Dw=w\frac{rV_r+V_u\cdot y}{V}-w\frac{DQ}{Q}.
\]
The complete coefficient of \(Dy\) in \(R[u]\) is
\begin{equation}\label{eq:acceleration-matrix}
 wg\left(I_k-\frac gQ yy^T\right),
 \qquad
 \left(I_k-\frac gQ yy^T\right)^{-1}=I_k+gyy^T.
\end{equation}
The inverse identity follows from \(Q=1+g|y|^2\). The scalars
\(g,Q\) are units and \(w\ne0\), so division in the Hahn field is
legitimate. We obtain a rational system
\[
 Du=y,\qquad Dy=\mathcal F(r,u,y),\qquad
 \mathcal F\in\R(r,u,y)^k.
\]
The field \(E=\R(r,u_1,\ldots,u_k,y_1,\ldots,y_k)\) is therefore
closed under \(D\) and has transcendence degree at most \(2k\) over
\(\R(r)\). Hence \(u_i,Du_i,\ldots,D^{2k}u_i\) satisfy a
nonzero polynomial relation over \(\R(r)\). Apply
Lemma~\ref{lem:scalar-support} to each nonzero component and take the
least common multiple of their denominators.
\end{proof}

This is a statement about formal support, not analytic convergence.
It also does not apply directly to a general analytic potential, whose
angular system need not be rational in its dependent variables.
Sections~\ref{sec:threshold}--\ref{sec:analytic-transfer} provide the
separate analytic argument.

\section{A fixed-denominator continuation criterion}
\label{sec:threshold}

All estimates in this section concern positive real \(r\) sufficiently close
to zero. Their constants may depend on the fixed angular dimension \(k\)
and on the current minimum certificate. No uniformity over an infinite
sequence of certificates is required. Write
\[
 K_N=\operatorname{Frac}\R\{r^{1/N}\}.
\]
Whenever a common denominator \(N\) is chosen, it includes the coefficients
of the current certificate and the monomial \(r^p\).

\begin{lemma}[Positive Puiseux factorization]\label{lem:factor}
Let \(M\in\operatorname{Mat}_k(K_N)\) be symmetric. Suppose that
\[
 M(r)\succeq c r^p I_k,\qquad c>0.
\]
There is an invertible matrix \(S\) over \(K_{2N}\) such that
\begin{equation}\label{eq:factor}
 M=S^TS,\qquad \|S^{-1}\|\le c^{-1/2}r^{-p/2}.
\end{equation}
\end{lemma}

\begin{proof}
Let \(\Delta_j\) be the leading principal minor of order \(j\), and set
\(\Delta_0=1\). Positive definiteness implies \(\Delta_j>0\). Symmetric
elimination, performed entirely in the field \(K_N\), gives
\[
 M=B^T\operatorname{diag}(d_1,\ldots,d_k)B,\qquad
 d_j=\Delta_j/\Delta_{j-1}>0,
\]
where \(B\) is invertible and triangular over \(K_N\).
Put \(t=r^{1/N}\). Each positive scalar germ has the form
\[
 d_j=t^{a_j}v_j(t),\qquad a_j\in\mathbb Z,\quad
 v_j\in\R\{t\},\quad v_j(0)>0.
\]
Consequently all its positive square roots belong to the same field
\(K_{2N}\). Set
\[
 S=\operatorname{diag}(\sqrt{d_1},\ldots,\sqrt{d_k})B.
\]
Then \(M=S^TS\), and
\(
 |Sh|^2=h^TMh\ge cr^p|h|^2
\)
proves the norm estimate.
\end{proof}

\begin{lemma}[A cubic estimate at a minimum]\label{lem:cubic-minimum}
Suppose \(F_r(0)=dF_r(0)=0\) and \(F_r(h)\ge0\) on a fixed ball about zero.
Assume that its third and fourth derivatives are \(O(r^m)\) there. Let
\[
 H=D_h^2F_r(0)\succeq0,\qquad T=D_h^3F_r(0).
\]
There are constants independent of small \(r\) such that
\begin{align}
 |T[h,h,h]|&\le C r^{m/2}|h|^2\sqrt{H[h,h]},
 \label{eq:cubic-bound}\\
 F_r(h)&\ge \tfrac14 H[h,h]-Cr^m|h|^4
 \label{eq:minimum-lower}
\end{align}
on a possibly smaller fixed ball.
\end{lemma}

\begin{proof}
For a unit vector \(e\), put \(A=H[e,e]\) and \(B=T[e,e,e]\).
Taylor's formula and nonnegativity imply, on a fixed two-sided interval,
\[
 0\le \frac{F_r(te)}{t^2}
 \le \frac A2+\frac B6t+C_0r^mt^2 .
\]
The quotient is interpreted continuously at \(t=0\). Since
\(|B|\le C r^m\) uniformly in \(e\), increasing the fixed constant \(C_0\)
places the vertex of the upper quadratic inside that interval.
Evaluation at the vertex gives \(B^2\le C r^m A\).
Homogeneity proves \eqref{eq:cubic-bound}. Applying Young's inequality to
the cubic term in Taylor's formula proves \eqref{eq:minimum-lower}.
\end{proof}

\begin{proposition}[Fixed-denominator continuation]\label{prop:threshold}
Let \(V(r,x)=U(r\phi(x))\), where \(U\) is real analytic of order \(m\).
Suppose a convergent Puiseux curve \(u(r)\) is a global minimum on \(K\) of
\[
 J(r,x)=V(r,x)+\frac{\lambda}{2}r^p|x|^2+\ell(r)\cdot x,
\]
where \(\ell\) is convergent Puiseux. Assume
\[
 \nu(u)=\tau>0,\quad V(r,u)=-\alpha r^p+o(r^p),\quad
 \lambda=C(d),\quad d\ge\tau,\quad
 f:=\nu (R[u])>p+d.
\]
If \(R[u]=0\), the curve is already an exact shape. Otherwise, if
\begin{equation}\label{eq:threshold}
 f>\frac{3p-m}{2},
\end{equation}
all subsequent matched corrections can be made with one fixed Puiseux
denominator. They either terminate at an exact convergent shape or
converge formally to an exact ordinary Puiseux shape with the same
leading curve and negative potential terms.
\end{proposition}

\begin{proof}
Fix \(c>0\), and choose a convergent Puiseux global minimum \(u+h\) of
\[
 J(r,x)+\frac c2r^p|x-u|^2+R[u]\cdot(x-u).
\]
Such selections exist by compact subanalytic minimization. The comparison
with \(x=u\), as in Proposition~\ref{prop:successor}, yields
\begin{equation}\label{eq:unweighted-displacement}
 |h|\le \frac{2|R[u]|}{cr^p}.
\end{equation}
In particular \(u+h\) is an interior minimum for small \(r\).

Set
\[
 F_r(h)=J(r,u+h)-J(r,u),\quad H=D_h^2F_r(0),\quad
 T=D_h^3F_r(0),\quad M=H+cr^pI_k.
\]
The center tends to the interior point zero of \(K\), so \(F_r\ge0\) on
a fixed ball. Its derivatives of order at least three are those of \(V\);
they satisfy the hypotheses of Lemma~\ref{lem:cubic-minimum}.
Choose \(N\) containing all current data. By Lemma~\ref{lem:factor},
\(M=S^TS\) over \(K_{2N}\). Introduce
\[
 \eta=f-p,\qquad \zeta=f-\frac p2,\qquad
 \gamma=f+\frac m2-\frac{3p}{2}>0,\qquad z=Sh=r^\zeta\xi.
\]
Notice that \(\eta>(p-m)/2\ge0\).

The critical-point equation is
\begin{equation}\label{eq:critical-correction}
 Mh+R[u]+\mathcal N(r,h)=0,
\end{equation}
where \(\mathcal N\) is the part of the gradient Taylor expansion of
degree at least two in \(h\). Left multiplication by \(S^{-T}\), followed
by division by \(r^\zeta\), normalizes the linear part to the identity.
Since \(H\preceq M\), \eqref{eq:cubic-bound} gives
\[
 |T[S^{-1}a,S^{-1}a,S^{-1}a]|
 \le C r^{m/2-p}|a|^3.
\]
Polarization of this symmetric trilinear form bounds all its mixed
coefficients with the same order. The normalized quadratic gradient
term thus has order at least
\(\zeta+m/2-p=\gamma\). A potential term of degree \(j\ge4\) contributes
order at least
\begin{equation}\label{eq:higher-normalized-order}
 m-\frac{jp}{2}+(j-2)\zeta
 =2\gamma+(j-4)\eta.
\end{equation}
We obtain
\begin{equation}\label{eq:normalized-implicit}
 \xi+b(r)+r^\gamma B_2(r)[\xi,\xi]
 +r^{2\gamma}B_{\ge3}(r,\xi)=0,
 \qquad b=r^{-\zeta}S^{-T}R[u].
\end{equation}
The vector \(b\) is bounded.

Here the remainder is jointly analytic after the single temporary
ramification; this requires more than coefficientwise order estimates.
To see it, put \(t=r^{1/(2N)}\) and
\[
 A(t)=r^{p/2}S^{-1}.
\]
The matrix \(A\) is a bounded convergent Laurent germ, hence is analytic
at \(t=0\). The fourth and higher order potential remainder has the form
\(r^m Q_4(t,h)\), where \(Q_4\) is jointly analytic and vanishes to order
at least four in \(h\). Substitution of
\(h=r^\eta A(t)\xi\) in \(\nabla_h Q_4\) extracts \(r^{3\eta}\), leaving
a jointly analytic function on each fixed bounded \(\xi\)-neighborhood.
The remaining normalization factor is
\[
 r^{m-p/2-\zeta+3\eta}=r^{2\gamma}.
\]
This proves joint analyticity of \(B_{\ge3}\). The cubic bound and the
removability of bounded Laurent poles prove analyticity of \(B_2\);
the same argument applies to \(b\). All exponents used here belong to
\((1/(2N))\mathbb Z\). Since \(\gamma>0\), the derivative with respect
to \(\xi\) at \(t=0,\xi=-b(0)\) in
\eqref{eq:normalized-implicit} is \(I_k\).

We must also identify the selected global minimum with this implicit
branch. Its minimum inequality is
\[
 F_r(h)+\frac c2r^p|h|^2+R[u]\cdot h\le0.
\]
By \eqref{eq:unweighted-displacement}, \(h=O(r^\eta)\).
The quartic error in \eqref{eq:minimum-lower} can be absorbed into
\((c/4)r^p|h|^2\), because
\[
 m+2\eta-p=2\gamma>0.
\]
It follows that
\[
 \frac14 h^TMh\le |R[u]\cdot h|
 \le |S^{-T}R[u]|\,|Sh|,
\]
and therefore
\begin{equation}\label{eq:weighted-minimum-bound}
 |Sh|\le4|S^{-T}R[u]|=O(r^\zeta).
\end{equation}
Every selected Puiseux minimum thus gives bounded \(\xi\). Equation
\eqref{eq:normalized-implicit} forces \(\xi\to-b(0)\); the analytic
implicit function theorem consequently identifies every such minimum
with its unique local branch.

The temporary denominator \(2N\) is not retained. In the original
coordinates iterate
\begin{equation}\label{eq:original-picard}
 h^{(0)}=-M^{-1}R[u],\qquad
 h^{(j+1)}=-M^{-1}R[u]-M^{-1}\mathcal N(r,h^{(j)}).
\end{equation}
Each iterate belongs to \(K_N^k\). Indeed the initial vector has positive
order; in normalized coordinates the iteration is the Picard map of
\eqref{eq:normalized-implicit}, so all iterates are bounded there and
have order at least \(\eta>0\) in the original coordinates. Analytic
substitution is consequently legitimate at every step. The normalized
map raises the valuation of differences by at least \(\gamma\), and
therefore converges formally to the implicit branch. Fixed
Laurent--Puiseux matrices and monomials change valuations by only finite
amounts, so the original iterates converge formally as well. Their limit
has no exponent outside \((1/N)\mathbb Z\). Since this limit is the
expansion of a branch convergent in \(r^{1/(2N)}\), that branch already
converges in the original parameter \(r^{1/N}\).

This argument holds for every fixed \(c>0\), and hence for the matched
constant selected in Proposition~\ref{prop:successor}. The new forcing,
center, and residual remain in \(K_N\). The successor inequalities are
\[
 d_{\mathrm{new}}\ge f-p,\qquad
 f_{\mathrm{new}}>p+d_{\mathrm{new}}\ge f,
\]
so \eqref{eq:threshold} persists. Unless an exact shape is reached, the
strictly increasing correction orders in the fixed discrete group
\((1/N)\mathbb Z\) tend to infinity; the residual orders do likewise.
The resulting formal limit solves the complete shape equation.
Analytic substitution, Euler differentiation, and division by a series
with fixed nonzero leading term are continuous in this ordinary formal
topology. The comparison estimates preserve the leading center and
potential terms. The proof also shows that the label \(d\) need not
equal \(\tau\), nor be the order of an earlier correction.
\end{proof}

\section{Transfer from polynomial to analytic potentials}
\label{sec:analytic-transfer}

The scalar differential-algebra argument for a polynomial potential
does not automatically extend to arbitrary analytic dependent-variable
coefficients. We instead extract a convergent minimum certificate from
one sufficiently high Taylor polynomial and transfer that certificate
back to the original potential.

\begin{lemma}[Transfer of a minimum certificate]\label{lem:minimum-transfer}
Let an objective on the fixed ball \(K\) be rational in its angular
variables, with coefficients that are convergent real Puiseux germs.
Assume its denominators are positive on \(K\). If a convergent Puiseux
curve is a global minimum for all sufficiently small positive real
\(r\), its Hahn expansion is a global minimum against all points of
\(K_{\HH_{\Q}}\).
\end{lemma}

\begin{proof}
The coefficients, the proposed minimum, and the monomial parameters
form a finite tuple. After clearing the positive denominators, global
minimality is a first-order formula in the language of ordered fields
with that tuple as parameters. Quantifier elimination for real closed
fields replaces this formula by a Boolean combination of finitely
many polynomial sign and equality conditions. Every resulting
convergent Puiseux expression has in \(\HH_{\Q}\) the same sign as its
germ on sufficiently small positive real \(r\); equality to zero
agrees as well. Hence the formula holds in the Hahn field.
\end{proof}

\begin{lemma}[Perturbing the potential at a quadratically pinned minimum]
\label{lem:perturbation}
Let \(V_1,V_2\) be analytic chart potentials such that, uniformly on \(K\),
\[
 V_2(r,x)-V_1(r,x)=O(r^\kappa).
\]
Suppose \(u\) is a convergent Puiseux global minimum of
\[
 J_1(r,x)=V_1(r,x)+\frac{\lambda}{2}r^p|x|^2+\ell(r)\cdot x,
\]
with
\[
 \nu(u)=\tau>0,\quad V_1(r,u)=-\alpha r^p+o(r^p),\quad
 \nu R_1[u]>p+\tau.
\]
Fix \(c>0\), and put
\[
 \lambda'=\lambda+c,\qquad \ell'=\ell-cr^pu.
\]
Let \(v\) be any convergent Puiseux global minimum of
\[
 J_2(r,x)=V_2(r,x)+\frac{\lambda'}2r^p|x|^2+\ell'(r)\cdot x.
\]
If \(\eta=(\kappa-p)/2>\tau\), then
\begin{align}
 \nu(v-u)&\ge\eta,\label{eq:potential-transfer-position}\\
 \nu\bigl(V_2(r,v)-V_1(r,u)\bigr)&\ge p+\tau+\eta,
 \label{eq:potential-transfer-energy}\\
 R_2[v]-R_1[u]&=O(r^{p+\eta}).
 \label{eq:potential-transfer-residual}
\end{align}
In particular \(v\) has the same leading curve term and negative
potential term as \(u\). Here \(R_i\) and \(\mathcal I_i\) are computed
with the potential \(V_i\).
\end{lemma}

\begin{proof}
The strengthened old objective
\(
 J'_1=V_1+\lambda'r^p|x|^2/2+\ell'\cdot x
\)
satisfies the global growth inequality
\[
 J'_1(r,x)-J'_1(r,u)\ge\frac c2r^p|x-u|^2.
\]
Write \(E=V_2-V_1\) and \(h=v-u\). Comparing the new minimum with \(u\)
gives
\[
 \frac c2r^p|h|^2\le E(r,u)-E(r,v)=O(r^\kappa),
\]
which proves \eqref{eq:potential-transfer-position}. The new minimum is
interior for small \(r\).

Set \(\Delta=J_2(r,v)-J'_1(r,u)\). The minimum comparison bounds it
above by \(E(r,u)\); the old minimum inequality bounds it below by
\(E(r,v)\). Thus \(\Delta=O(r^\kappa)\). Stationarity of \(J'_1\) at \(u\)
and exact expansion of the common regularizer give
\begin{equation}\label{eq:potential-transfer-exact}
 V_2(r,v)-V_1(r,u)
 =\Delta+\partial_uV_1(r,u)\cdot h
 -\frac{\lambda'}2r^p|h|^2.
\end{equation}
The residual assumption and the leading inertia estimate imply
\(\partial_uV_1(r,u)=O(r^{p+\tau})\).
Since \(\kappa=p+2\eta\) and \(\eta>\tau\),
\eqref{eq:potential-transfer-exact} proves
\eqref{eq:potential-transfer-energy}.

For completeness, the inertia comparison remains valid when the
potential changes. The denominators
\(Q_i=1+g(u_i)|Du_i|^2\), with \(u_1=u,u_2=v\), satisfy
\(Q_2-Q_1=O(r^{\tau+\eta})\).
Consequently the corresponding energy factors obey
\[
 w_2-w_1=O(r^{p+\tau+\eta}).
\]
Euler differentiation preserves these Puiseux order bounds. In the
inertia difference, terms linear in \(h,Dh,D^2h\) have order at least
\(p+\eta\); terms involving changes of \(w,Dw\), the metric, or the
connection have additional positive-order factors. Therefore
\[
 \mathcal I_2[v]-\mathcal I_1[u]=O(r^{p+\eta}).
\]
The two exact critical-point equations use the same \(\lambda',\ell'\),
and hence
\[
 \partial_uV_2(r,v)-\partial_uV_1(r,u)=-\lambda'r^ph.
\]
Adding the inertia difference proves
\eqref{eq:potential-transfer-residual}. This argument does not estimate
the force perturbation merely from the size of the potential
perturbation; the two stationary equations supply the required
identity.
\end{proof}

\begin{proposition}[Analytic transfer]\label{prop:analytic}
Under the hypotheses and notation of the spherical minimum
construction, an arbitrary real analytic potential admits a complete
ordinary real formal Puiseux shape with negative leading potential.
\end{proposition}

\begin{proof}
If Proposition~\ref{prop:initial} gives a radial branch, or an exact
convergent shape, the assertion follows immediately. Otherwise let
\(u_0\) be its convergent matched center, of order \(\tau>0\), with
\[
 \lambda_0=C(\tau),\qquad f_0=\nu R_U[u_0]>p+\tau.
\]
If \(f_0>(3p-m)/2\), apply Proposition~\ref{prop:threshold}.
It remains to treat
\begin{equation}\label{eq:hard-initial-range}
 p+\tau<f_0\le p+b,\qquad
 b=\frac{p-m}{2}.
\end{equation}
In particular \(b>0\) and \(0<\tau<b\).
Put \(B=3b+1\), and choose the Taylor polynomial \(P\) of \(U\)
through an integer degree \(J\) such that
\begin{equation}\label{eq:taylor-choice}
 J+1>p+2B=4p-3m+2.
\end{equation}
Let
\[
 V_P(r,x)=P(r\phi(x)),\qquad H_*=\frac{J+1-p}{2}>B.
\]
Uniformly on the sphere, and hence on \(K\), the difference of the
potentials is \(O(r^{J+1})\). Since \(J+1>p\ge m\), the polynomial
has the same \(m,p,\alpha\) as \(U\).

\medskip\noindent
\emph{Initialization with a bounded label.}
Choose
\[
 \tau<\bar d<f_0-p\le b.
\]
Apply Lemma~\ref{lem:perturbation} to the original certificate at \(u_0\),
with
\[
 V_1=V_U,\qquad V_2=V_P,\qquad
 c=C(\bar d)-\lambda_0>0,\qquad \kappa=J+1.
\]
The strengthening in that lemma raises the coefficient exactly to
\(C(\bar d)\) and adds a positive quadratic growth margin.
The new actual Puiseux minimum \(u_P\)
has order \(\tau\), the same negative leading potential, and residual
order exactly \(f_0\), because
\[
 p+H_*>p+B>p+b\ge f_0.
\]
It is a matched certificate with label \(\bar d<b\).
Its forcing is a convergent Puiseux germ, although not necessarily
algebraic over \(\R(r)\). Lemma~\ref{lem:minimum-transfer} nevertheless
makes it a global minimum against all Hahn test points. We may
therefore begin the polynomial construction of
Proposition~\ref{prop:hahn} from this certificate.

\medskip\noindent
\emph{The first threshold stage.}
Run that polynomial construction to a full solution
\(\widehat u_P\). By Proposition~\ref{prop:puiseux},
\(\widehat u_P\) is an ordinary formal Puiseux vector.
Let \(\beta_*\) be the first stage at which the residual vanishes or
\[
 f_{\beta_*}>p+b.
\]
Such a stage exists. Indeed, before any such stage, \(f_\beta\le p+b\)
and \(f_\beta>p+d_\beta\) imply \(d_\beta<b\), so an unbounded-label
completion cannot occur without an earlier crossing.

The first crossing has the uniform bound
\begin{equation}\label{eq:first-crossing-label}
 d_{\beta_*}\le3b.
\end{equation}
For a successor, Proposition~\ref{prop:successor} gives
\(
 d_{\beta_*}\le f_{\mathrm{prev}}-m\le p+b-m=3b.
\)
For a bounded-limit stage, the preceding labels are all less than \(b\),
so their supremum is at most \(b\). The initial label is already less
than \(b\).

If the residual at this stage is zero, take
\(\widehat u_P=u_{\beta_*}\) and write \(f_{\beta_*}=+\infty\).
Otherwise, the all-pair estimates of Lemma~\ref{lem:coherence}, passed
to the final solution if necessary, give
\begin{align}
 \nu(\widehat u_P-u_{\beta_*})&\ge f_{\beta_*}-p,\label{eq:final-anchor-position}\\
 \nu\bigl(V_P(r,\widehat u_P)-V_P(r,u_{\beta_*})\bigr)
 &\ge f_{\beta_*}+\tau,\label{eq:final-anchor-energy}\\
 \nu\bigl(\mathcal I_P[\widehat u_P]-\mathcal I_P[u_{\beta_*}]\bigr)
 &\ge f_{\beta_*}.\label{eq:final-anchor-inertia}
\end{align}
We use the full polynomial solution to establish the next finite
truncation. Proposition~\ref{prop:threshold} is not being applied to a
possibly non-Puiseux intermediate Hahn center.

\medskip\noindent
\emph{Extraction of a convergent certificate.}
Choose a real number \(d_1\) such that
\begin{equation}\label{eq:first-new-label}
 \max(d_{\beta_*},b)<d_1<\min(f_{\beta_*}-p,B),
 \qquad \Lambda_1=C(d_1).
\end{equation}
This interval is nonempty by \eqref{eq:first-crossing-label},
\(3b<B\), and \(f_{\beta_*}-p>\max(d_{\beta_*},b)\).
Define
\[
 b_{\beta_*}=\ell_{\beta_*}-(\Lambda_1-\lambda_{\beta_*})r^pu_{\beta_*}.
\]
Stationarity of the old certificate gives the useful exact expression
\begin{equation}\label{eq:anchor-forcing}
 b_{\beta_*}=\mathcal I_P[u_{\beta_*}]-\Lambda_1r^pu_{\beta_*}-R_P[u_{\beta_*}].
\end{equation}
By \eqref{eq:final-anchor-position}--\eqref{eq:final-anchor-inertia},
the coefficients of \(b_{\beta_*}\) below \(f_{\beta_*}\) agree with those of
\[
 \mathcal I_P[\widehat u_P]-\Lambda_1r^p\widehat u_P.
\]
The latter is an ordinary formal Puiseux vector. Choose a rational
number \(T\) with
\begin{equation}\label{eq:forcing-cutoff}
 p+d_1<T<\min(f_{\beta_*},p+B).
\end{equation}
The strict truncation
\(\widetilde b=(b_{\beta_*})_{<T}\) is therefore a finite Puiseux
polynomial. This finiteness is obtained from the completed solution,
not from a support assertion for an incomplete differential equation.

Choose an actual convergent Puiseux global minimum \(v\) of
\[
 \widetilde J(r,x)
 =V_P(r,x)+\frac{\Lambda_1}{2}r^p|x|^2+\widetilde b(r)\cdot x
\]
on \(K\). Compact subanalytic selection gives such a curve; by
Lemma~\ref{lem:minimum-transfer} it is also a Hahn global minimum.
The old center \(u_{\beta_*}\) minimizes
\[
 J_{\beta_*}^+(r,x)
 =V_P(r,x)+\frac{\Lambda_1}{2}r^p|x|^2+b_{\beta_*}(r)\cdot x
\]
with the global growth margin
\[
 J_{\beta_*}^+(r,x)-J_{\beta_*}^+(r,u_{\beta_*})
 \ge\frac{\Lambda_1-\lambda_{\beta_*}}{2}r^p|x-u_{\beta_*}|^2.
\]
Write \(h=v-u_{\beta_*}\) and
\(\Delta b=\widetilde b-b_{\beta_*}\). Since \(\nu(\Delta b)\ge T\),
comparison of these two minima yields
\begin{equation}\label{eq:extracted-position}
 |h|\le\frac{2|\Delta b|}
 {(\Lambda_1-\lambda_{\beta_*})r^p},
 \qquad \nu(h)\ge T-p.
\end{equation}
Thus \(v\) is interior and has the same leading center term.

To control its residual, retain the unstrengthened old objective
\(J_{\beta_*}\). The same comparison gives
\[
 0\le J_{\beta_*}(r,v)-J_{\beta_*}(r,u_{\beta_*})
 \le|\Delta b|\,|h|.
\]
Exact expansion of its regularizer gives
\[
 \begin{split}
 V_P(r,v)-V_P(r,u_{\beta_*})
 ={}&J_{\beta_*}(r,v)-J_{\beta_*}(r,u_{\beta_*})\\
 &+\partial_uV_P(r,u_{\beta_*})\cdot h
 -\frac{\lambda_{\beta_*}}{2}r^p|h|^2 .
 \end{split}
\]
Here \(\partial_uV_P(r,u_{\beta_*})=O(r^{p+\tau})\), while
\(T-p>d_1>\tau\). Therefore
\[
 \nu\bigl(V_P(r,v)-V_P(r,u_{\beta_*})\bigr)\ge T+\tau,
 \qquad
 \nu\bigl(\mathcal I_P[v]-\mathcal I_P[u_{\beta_*}]\bigr)\ge T,
\]
the second estimate following from Lemma~\ref{lem:inertia}.
Subtracting the two strengthened critical-point equations gives
\[
 R_P[v]=R_P[u_{\beta_*}]
 +\mathcal I_P[v]-\mathcal I_P[u_{\beta_*}]
 -\Lambda_1r^ph-\Delta b.
\]
Every term on the right has order at least \(T\). Hence
\begin{equation}\label{eq:extracted-residual}
 \nu R_P[v]\ge T>p+d_1>p+b.
\end{equation}
We have obtained a convergent minimum certificate; no convergence of
the full formal solution \(\widehat u_P\) has been assumed.

\medskip\noindent
\emph{Return to the original analytic potential.}
Choose
\[
 d_1<d_2<T-p<B,\qquad \Lambda_2=C(d_2).
\]
Increase the quadratic coefficient from \(\Lambda_1\) to \(\Lambda_2\)
and replace \(\widetilde b\) by
\[
 \ell_2=\widetilde b-(\Lambda_2-\Lambda_1)r^pv.
\]
The same center \(v\) remains a global minimum for \(P\), now with a
positive quadratic growth margin. This is precisely the strengthening
in Lemma~\ref{lem:perturbation}: apply that lemma to the certificate
with coefficient \(\Lambda_1\) and forcing \(\widetilde b\), taking
\(c=\Lambda_2-\Lambda_1\), \(V_1=V_P\), \(V_2=V_U\), and
\(\kappa=J+1\). It supplies a convergent Puiseux minimum \(z\)
for the original analytic potential, with coefficient \(\Lambda_2\)
and forcing \(\ell_2\), such that
\[
 \nu(z-v)\ge H_*>B>d_2,\qquad
 R_U[z]=R_P[v]+O(r^{p+H_*}).
\]
It preserves the center order \(\tau\) and the negative potential
term \(-\alpha r^p\). By \eqref{eq:extracted-residual},
\[
 \nu R_U[z]\ge\min(T,p+H_*)>p+d_2>p+b
 =\frac{3p-m}{2}.
\]
If this residual is zero, \(z\) is already an exact shape. Otherwise
it satisfies every hypothesis of Proposition~\ref{prop:threshold},
with label \(d_2\ge\tau\), and that proposition produces the required
complete ordinary formal Puiseux shape for \(U\).
\end{proof}

\section{Realization in the original phase space}\label{sec:realization}

\subsection{Formal invariant curves and the external theorem}
For a nonzero formal parametrization
\(\Gamma(s)\in(s\R[[s]])^d\), set
\begin{equation}\label{eq:formal-curve-ideal}
 I_{\mathcal C}=\ker\bigl(\R[[z_1,\ldots,z_d]]\longrightarrow
 \R[[s]],\ F\longmapsto F(\Gamma(s))\bigr).
\end{equation}
Substitution is well defined coefficientwise because all coordinates
have positive order. The kernel is prime. Invariance means
\(\widehat X(I_{\mathcal C})\subset I_{\mathcal C}\), where
\(\widehat X\) is the Taylor derivation of the vector field. In our
application we will prove the stronger parametric identity
\(\widehat X\circ\Gamma=a(s)\Gamma'(s)\), as well as dimension one
of the quotient ring and \(\widehat X\circ\Gamma\ne0\).

A real formal half-branch specifies a side of a primitive parameter,
equivalently its oriented iterated tangents under successive spherical
point blow-ups. An actual curve is asymptotic to it when their oriented
iterated tangents agree. In nonsingular graph coordinates
\((x,\widehat h(x))\), on the positive side, this is
\begin{equation}\label{eq:formal-graph-asymptotics}
 h(x)-j_N\widehat h(x)=O(x^{N+1})\quad(x\downarrow0)
 \qquad\text{for every }N,
\end{equation}
as in \cite[Section~5, (5.1)]{LeGalSanz2026}. Each comparison uses a
finite polynomial; convergence of \(\widehat h\) is not required.
The constants and neighborhoods may depend on \(N\), whereas the
actual curve is the same for all \(N\).

For the following statement, a horn neighborhood has opening \(k\)
if \(k\) is the least number of spherical blow-ups along the half-branch
for which the closure of its inverse image is a neighborhood of the
selected lifted point \cite[Definition~5.1]{LeGalSanz2026}.

\begin{externaltheorem}[Le Gal--Sanz S\'anchez, {\cite[Theorem~5.1]{LeGalSanz2026}}]
Let \(X\) be \(C^\infty\) near \(a\in\R^d\), and let
\(\mathcal C\) be an invariant irreducible real formal curve at
\(a\), not contained in the formal singular locus. Fix a half-branch
\(\mathcal C^\varepsilon\). For every neighborhood \(U\) of \(a\)
there is \(k_0\) such that every \(k\ge k_0\) admits a horn
\(V^\varepsilon\subset U\setminus\operatorname{Sing}(X)\) of
opening \(k\), and a connected, relatively closed, positive-dimensional
\(C^0\) submanifold \(S^\varepsilon\subset V^\varepsilon\), locally
invariant under \(X\). Trajectories through \(S^\varepsilon\)
accumulate at \(a\) and are asymptotic to \(\mathcal C^\varepsilon\).
Trajectories through \(V^\varepsilon\setminus S^\varepsilon\)
leave \(V^\varepsilon\) in finite positive and negative time.
\end{externaltheorem}

In the last clause we write the trajectory through the quantified point
of \(V^\varepsilon\setminus S^\varepsilon\); the printed statement
has ``through \(a\)'', an apparent typographical slip clarified by
Step~4 of its proof. That clause is not used here.

\begin{theorem}[Realization consequence]\label{thm:realization-input}
Let \(X\) be a \(C^\infty\) vector field near zero in a
finite-dimensional real space. An invariant irreducible real formal
curve through zero, not contained in the formal singular locus, has a
nonconstant trajectory germ tending to zero and asymptotic to any
specified real half-branch.
\end{theorem}
\begin{proof}
Apply the external theorem with \(a=0\). Its positive-dimensional
set \(S^\varepsilon\) is nonempty. The trajectory through any of
its points has the required asymptotics and is nonconstant because
\(V^\varepsilon\) avoids the singular set. Restrict to its germ at
the limiting endpoint.
\end{proof}

The consequence fixes a geometric half-branch, not a direction of
physical time. We determine that direction after constructing the
formal curve and checking every hypothesis for the original field.

\subsection{Proof of the main theorem}
\begin{proof}[Proof of Theorem~\ref{thm:main}]
The radial alternatives, including \(n=1\), were constructed in
Section~\ref{subsec:radial}. In the remaining case,
Proposition~\ref{prop:analytic} supplies an ordinary formal Puiseux
shape \(\widehat u\) satisfying
\[
 R[\widehat u]=0,\qquad V(r,\widehat u)=-\alpha r^p+o(r^p),
 \qquad p\ge2.
\]
Let \(w=w[\widehat u]\), now a formal series with leading term
\(2\alpha r^p\). Enlarge the common denominator \(N\), by at most
a factor two, to include its square root with positive leading
coefficient. Define
\begin{align*}
 q(r)&=r\phi(\widehat u(r)),& v(r)&=-q_r(r)\sqrt{w(r)},\\
 \Gamma(s)&=(q(s^N),v(s^N)),&
 a(s)&=-\frac{\sqrt{w(s^N)}}{Ns^{N-1}}.
\end{align*}
All coordinates of \(\Gamma\) have positive order, since
\(q_r=e_n+o(1)\), and
\(\ord_s a=N(p/2-1)+1\ge1\). Thus both \(\Gamma\) and \(a\)
are ordinary formal power series. In particular,
\begin{equation}\label{eq:phase-leading}
 q(s^N)=s^Ne_n+o(s^N),\qquad
 v(s^N)=-\sqrt{2\alpha}\,s^{Np/2}e_n+o(s^{Np/2}).
\end{equation}

\emph{Recovery of the radial equation.}
Set
\[
 Z=wq_{rr}+\tfrac12w_rq_r+\nabla U(q).
\]
The angular Euler--Lagrange identity is
\(R_a[\widehat u]=r\,\partial_a\phi\cdot Z\). Thus \(Z\)
is orthogonal to all angular tangents. Their Gram matrix is \(gI_k\),
so \(Z=\beta\phi\) over the formal fraction field. On the other
hand, \(|q_r|^2=Q[\widehat u]\) and \(wQ=-2V\), giving
\[
 \tfrac12w|q_r|^2+U(q)=0.
\]
Differentiation in \(r\) gives \(Z\cdot q_r=0\). Since
\(\phi\cdot q_r=1\), it follows that \(\beta=0\). Consequently,
for the original Newton field \(X_U(q,v)=(v,-\nabla U(q))\),
\begin{equation}\label{eq:formal-invariance}
 X_U\circ\Gamma=a(s)\Gamma'(s).
\end{equation}
Indeed \(a,dr/ds=-\sqrt w\), and
\((-\sqrt w)d(-q_r\sqrt w)/dr=wq_{rr}+w_rq_r/2\).
No time-rescaled field is substituted into the realization theorem.

\emph{Verification of the formal-curve hypotheses.}
The field \(X_U\) is analytic on a neighborhood in \(\R^{2n}\)
and vanishes at \((0,0)\). Define \(I_{\mathcal C}\) by
\eqref{eq:formal-curve-ideal} in these phase coordinates. Its
substitution image \(B\) contains \(A=\R[[Q_n(s)]]\), where
\(Q_n(s)=q_n(s^N)=s^N\eta(s)\) and \(\eta(0)=1\).
The formal parameter change \(\sigma=s\eta(s)^{1/N}\) gives
\(Q_n=\sigma^N\). Hence \(\R[[s]]=\R[[\sigma]]\) is a free
\(A\)-module of rank \(N\). Its \(A\)-submodule \(B\) is finite
because \(A\) is Noetherian. Therefore \(\dim B=1\).
The prime kernel thus defines an irreducible real formal curve,
independently of possible multiple covering by \(\Gamma\).

For \(F\in I_{\mathcal C}\), the formal chain rule gives
\begin{equation}\label{eq:newton-ideal-invariance-proof}
 (\widehat X_UF)(\Gamma(s))
 =dF_{\Gamma(s)}\bigl(a(s)\Gamma'(s)\bigr)
 =a(s)\frac{d}{ds}F(\Gamma(s))=0.
\end{equation}
Thus \(\widehat X_U(I_{\mathcal C})\subset I_{\mathcal C}\).
The first block of \(\widehat X_U\circ\Gamma\) is \(v(s^N)\),
whose \(n\)-th component has nonzero leading coefficient in
\eqref{eq:phase-leading}. The curve is therefore not contained in
the formal singular locus. Finally choose the real half-branch
specified by \(s>0\); a positive multiple covering has the same
oriented iterated tangents. These facts verify smoothness, real
irreducibility, invariance, nonvanishing, and the chosen side required
by Theorem~\ref{thm:realization-input}.

\emph{Actual trajectory and physical time.}
The realization consequence gives a nonconstant actual trajectory of
\(X_U\) tending to zero along the chosen half-branch. Its finite
leading data imply \(q_n>0\) and \(v_n<0\) on a sufficiently short
representative. To justify using signs, resolve the formal curve by
finitely many point blow-ups and use graph coordinates there. A
coordinate with nonzero formal leading term has its sign fixed by a
sufficiently long finite jet in \eqref{eq:formal-graph-asymptotics};
blowing down preserves its sign as a function on the original curve.
This uses no numerical evaluation of a divergent series.

The exact equation \(\dot q_n=v_n<0\) shows that the endpoint
\(q_n\to0\) is the future endpoint: a decreasing positive function
cannot tend to zero at the past endpoint. It cannot reach zero at a
finite physical time. Otherwise extending the limiting state by the
constant zero state gives a solution through that time, by the integral
form of \(\dot z=X_U(z)\). Local uniqueness for the analytic field,
also backwards in time, would make the approaching trajectory constant.
Thus the limit occurs as \(t\to+\infty\).

Energy is constant on this actual orbit and has limit zero, hence
vanishes identically. Since \(v_n<0\), its velocity is nonzero and
\(U(q)=-|v|^2/2<0\). Taking a sufficiently late point as time zero
gives all assertions in \eqref{eq:mainorbit} on \([0,+\infty)\).

\emph{Forward instability.}
This step applies also to the radial solutions. Choose \(T_j\to+\infty\)
and set
\begin{equation}\label{eq:forward-escape}
 Q_j(t)=q(T_j-t),\qquad 0\le t\le T_j.
\end{equation}
Then \(Q_j\) solves the same Newton equation, and
\[
 (Q_j(0),\dot Q_j(0))=(q(T_j),-\dot q(T_j))\longrightarrow0,
 \quad
 (Q_j(T_j),\dot Q_j(T_j))=(q(0),-\dot q(0))=:z_*\ne0.
\]
Take \(\epsilon_0=|z_*|/2\). For every \(\delta>0\), a sufficiently
large \(j\) supplies \eqref{eq:forward-instability}.
\end{proof}

\section{Curved valleys and the continuation threshold}\label{sec:examples}

Consider
\begin{equation}\label{eq:valley-example}
 U_p(x,y,z)=(x-y^2)^2-y^p+z^2,\qquad p\in\{6,8,10\}.
\end{equation}
The origin is an isolated complex critical point: the \(x\)- and
\(z\)-equations give \(x=y^2\), \(z=0\), after which the
\(y\)-equation gives \(y=0\). No fixed ray is negative near zero.
Indeed, along \(r(a,b,c)\) the leading term is
\(r^2(a^2+c^2)\) unless \(a=c=0\), in which case the potential
is \(r^4-r^p>0\). The curved arc \((s^2,s,0)\), however, has
potential \(-s^p\).

Here \(m=2\), and the spherical minimum has leading term \(-r^p\).
The lower bound is \(U_p\ge-|y|^p\ge-r^p\). For the upper bound,
use the point \((y^2,y,0)\) on the sphere, where
\(y=r+O(r^3)\). Center the stereographic chart at the positive
\(y\)-axis. On the plane \(z=0\),
\[
 x=\frac{ru}{1+u^2/4},\qquad
 y=\frac{r(1-u^2/4)}{1+u^2/4}.
\]
The Hessian of the leading angular potential, divided by \(r^2\),
is \(2I\) at zero and remains positive on a small fixed ball.
Thus the regularized objective has a unique convergent interior
minimum; symmetry sets the other angular coordinate to zero.

The initial matched center has \(\tau=1\),
\(\lambda=C(1)=p+4\), and
\[
 u(r)=r-\frac34r^3+O(r^5).
\]
For example, substituting \(u=r+cr^3+O(r^5)\) gives
\(x-y^2=(c+3/4)r^4+O(r^6)\). The order-five stationary coefficient
is \(2(c+3/4)\), whereas the negative potential term and regularizer
first contribute at order \(p+1\ge7\), proving \(c=-3/4\).
The exact identity \(V_u=-(p+4)r^pu\) also gives
\(x-y^2=O(r^p)\). Therefore
\begin{align*}
 V(r,u)&=-r^p+\frac p2r^{p+2}+O(r^{p+4}),&
 g(u)&=1-\frac12r^2+O(r^4),\\
 Q[u]&=1+r^2+O(r^4),&
 w[u]&=2r^p-(p+2)r^{p+2}+O(r^{p+4}),\\
 Du&=r-\frac94r^3+O(r^5),&
 D^2u&=r-\frac{27}{4}r^3+O(r^5).
\end{align*}
Using \(\Gamma_{\rm cov}(u)[Du,Du]=-\rho^{-3}u(Du)^2/2\)
in \eqref{eq:inertia} gives
\[
 R[u]=-\frac{p^2+12p+48}{2}\,r^{p+3}+O(r^{p+5}).
\]
In particular,
\[
\begin{array}{c|c|c|c}
 p&\text{leading residual}&f=p+3&(3p-m)/2\\ \hline
 6&-78r^9&9&8\\
 8&-104r^{11}&11&11\\
 10&-134r^{13}&13&14
\end{array}
\]
Only the first row immediately satisfies the strict threshold.
For all three, one further correction has displacement order
\(\sigma=f-2=p+1\), since the angular Hessian starts with \(2r^2\).
Its leading coefficient is \((p^2+12p+48)/4\), independently of the
fixed positive regularization increment. Matching the new coefficient
to \(C(p+1)=(p+1)(3p+4)\) cancels order \(2p+1\), so the new
residual lies above the threshold in each case.

These examples test the finite correction and the strictness of the
threshold, not bounded ordinal accumulation. They also are not
intrinsically non-quasihomogeneous function singularities:
\((X,Y,Z)=(x-y^2,y,z)\) changes the potential to
\(X^2+Z^2-Y^p\). It changes the kinetic metric, so the preceding
Euclidean calculation is not reduced to a radial Newton problem by
that substitution.

\section{An obstruction to exact Euler-field certificates}\label{sec:euler}

The trajectory construction does not require a field \(A\) with
\(A(U)=U\). Such an identity can fail even under complex isolation
and under a weak relative differentiability requirement. On the
nonopen set \(X=\{U<0\}\cup\{0\}\), we impose only
\begin{equation}\label{eq:relative-expansion}
 A(x)=a+Lx+o(|x|)\qquad(x\to0,\ x\in X).
\end{equation}
Every restriction of a neighborhood \(C^1\) field has this property;
no converse extension assertion is made.

\begin{proposition}\label{prop:euler}
For every \(n\ge2\), the polynomial
\begin{equation}\label{eq:euler-counterexample}
 U_n(x)=-x_1^4-x_1^2x_2^2-x_2^6+\sum_{j=3}^n x_j^2
\end{equation}
has an isolated complex critical point at zero and is not a real local
minimum there. No field on \(X_n=\{U_n<0\}\cup\{0\}\) satisfying
\(A(U_n)=U_n\) admits \eqref{eq:relative-expansion}.
\end{proposition}
\begin{proof}
The complex gradient equations give \(x_j=0\) for \(j\ge3\), and
\[
 x_1(2x_1^2+x_2^2)=0,\qquad x_2(x_1^2+3x_2^4)=0.
\]
If either coordinate is zero, both are zero. Otherwise
\(x_2^2=1/6\), which excludes critical points sufficiently close to
zero. The real punctured plane \(x_3=\cdots=x_n=0\) lies in
\(\{U_n<0\}\).

Suppose \eqref{eq:relative-expansion} holds. On the first axis the
identity forces \(A_1=x_1/4\), and on the second it forces
\(A_2=x_2/6\). Hence
\[
 a_1=a_2=0,\qquad L_{11}=\frac14,\qquad L_{22}=\frac16.
\]
On the negative plane the other gradient components vanish. The
degree-four part contributed by the linear expansion is
\[
 -4L_{11}x_1^4-(4L_{12}+2L_{21})x_1^3x_2
 -2(L_{11}+L_{22})x_1^2x_2^2-2L_{12}x_1x_2^3.
\]
The remainder is \(o(|x|^4)\), because the gradient begins in degree
three. Fix a nonzero \(\xi\in\R^2\), substitute \(x=s\xi\),
divide by \(s^4\), and let \(s\downarrow0\). Equality for every
\(\xi\) gives equality of quartic polynomials. The mixed coefficient
requires \(2L_{11}+2L_{22}=1\), contradicting
\(1/2+1/3=5/6\).
\end{proof}

No condition on \(\operatorname{Sym}L\) was used. Thus the example
also excludes exact certificates with positive semidefinite linear
part, including full-neighborhood \(C^1\) and analytic certificates.
In the analytic category, an exact identity on a neighborhood implies
membership in the Jacobian ideal after complexification. Under complex
isolation, Saito's criterion \cite{Saito1971} relates this membership
to analytic equivalence to a weighted-homogeneous germ. The proposition
is stronger in regularity for this particular example because it only
uses relative first-order data on the negative set.

This does not obstruct Palamodov's weaker criterion
\cite[Theorem~2]{Palamodov2020}. In fact, set
\(B=(x_1/4,x_2/6,0,\ldots,0)\),
\(P=x_1^4+x_1^2x_2^2+x_2^6\), and
\(S=\sum_{j=3}^n x_j^2\). On \(U_n=-P+S<0\),
\[
 \frac{B(U_n)}{U_n}
 =\frac{x_1^4+\frac56x_1^2x_2^2+x_2^6}{P-S}
 \ge\frac56,
 \qquad
 \operatorname{Sym}DB(0)=\diag(1/4,1/6,0,\ldots,0)\succeq0.
\]
In dimension one, by contrast, if \(U(x)=x^p a(x)\), \(a(0)\ne0\),
then
\[
 A(x)=\frac{x a(x)}{p a(x)+x a'(x)}
\]
is analytic near zero, satisfies \(A(U)=U\), and has
\(A'(0)=1/p>0\).

\bibliographystyle{amsplain}
\bibliography{references_submission}
\end{document}